\documentclass[12pt]{amsart}
\usepackage{amscd}
\usepackage{verbatim}
\usepackage{amssymb, amsmath, amsthm, amscd,ifthen}
\usepackage[dvips]{graphics}
\usepackage[cp866]{inputenc}
\usepackage{graphicx}
\usepackage{epsfig}
\usepackage{tikz}
\usepackage{color}

\usepackage{amsmath,amssymb,amscd,}
\usepackage[all]{xy}

\usepackage[dvips]{graphics}
\usepackage[cp866]{inputenc}
\usepackage{graphicx}
\usepackage{epsfig}
\usepackage[hang]{subfigure}

\theoremstyle{plain}
\newtheorem{theorem}{Theorem}
\newtheorem{lemma}{Lemma}

\newtheorem{proposition}{Proposition}

\theoremstyle{definition}

\theoremstyle{plain}
\newtoks\thehProclaim
\newtheorem*{Proclaim}{\the\thehProclaim}

\theoremstyle{definition}
\newtoks{\thehRemark}
\newtheorem*{Remark}{\the\thehRemark}

\renewcommand{\leq}{\leqslant}
\renewcommand{\geq}{\geqslant}

\begin{document}

\title{Area-perimeter duality in polygon spaces}

\author{Giorgi Khimshiashvili, Gaiane Panina, Dirk Siersma}

\address{G. Khimshiashvili: Ilia State University, Tbilisi; gogikhim@yahoo.com; G. Panina: St. Petersburg Department of Steklov Mathematical Institute; St. Petersburg State University;  gaiane-panina@rambler.ru;  D. Siersma: Utrecht University, Department of Mathematics; d.siersma@uu.nl}



\begin{abstract}
{Two natural foliations, guided by area and perimeter, of the configurations spaces of planar polygons are considered
and the topology of their leaves is investigated in some detail. In particular, the homology
groups and the homotopy type of leaves  are determined. The homology groups of the spaces of polygons
with fixed area and perimeter are also determined. Besides, we extend
the  classical isoperimetric duality to all critical points. In
conclusion a few general remarks on dual extremal problems in polygon spaces and beyond are
given. }

\end{abstract}

\maketitle

\section{Introduction}

The aim of the present paper is to explicate  and generalize the results of the classical isoperimetric problem in polygon spaces given in \cite{Al, ropepaper}. In particular, we complement results of \cite{ropepaper} on the critical points of area on the space of polygons with fixed perimeter and generalize the results of the structure of area foliation in the space of triangles given in \cite{Giorg}.
One of the essential novel ingredients of our approach is a conceptual generalization of the well-known duality between area and perimeter in the classical isoperimetric problem { (see e.g. \cite{Leichtweiss})}. Namely, we extend the idea of duality to all critical points of the two dual extremal problems in question. Technically, this involves comparison of the Morse indices of two smooth functions at a generic tangent point of their level surfaces and applying a type of argument in the spirit of 'rabattement de diagramme de Cerf' widely used in singularity theory.

This, in particular, enables us to show that the perimeter function on the space of polygons with constant non-zero area is { level-preserving homeomorphic}  to the area function on the space of polygon with constant perimeter and positive area  (Theorem \ref{ThmRotate}), which may be considered as one of the main results of the present paper. We describe  also topological structure of these spaces.

 Another new result (which is crucial for the whole paper) is determining the homology groups of the spaces of $n$-gons with both area and perimeter fixed (Theorem \ref{ThmOdd}).

\section{Preliminaries}
An $n$\textit{-gon} is an $n$-tuple of points  $(p_1,...,p_n)\in (\mathbb{R}^2)^n$  modulo diagonal action of orientation preserving isometries of the plane. Some  (but not all) of $p_i$ may coincide.

The space $\mathcal{C}=\mathcal{C}^n$ { of all $n$-gons for a fixed $n$} is equipped by two natural functions
$$\mathcal{A},\mathcal{P}:\mathcal{C}\rightrightarrows \mathbb{R}:$$
\begin{itemize}
  \item The \textit{perimeter} of a polygon is (as usual)  $$\mathcal{P}(p_1,...,p_n)=|p_1p_2|+|p_2p_3|+...+|p_np_1|.$$
  \item The \textit{oriented area} is  defined as
$$2\mathcal{A}(p_1,...,p_n)=x_1y_2-x_2y_1+...+x_ny_1-x_1y_n,$$

where $p_i=(x_i,y_i)$.
\end{itemize}

One sees immediately that perimeter of a polygon is always positive, whereas the area can be {also} zero or negative.
In particular,  the set $\{\mathcal{A}=0\}$ splits $\mathcal{C}$ into two mutually symmetric parts.

\medskip

Oriented area as a Morse function has been studied in various settings:  for configuration spaces of flexible polygons (those with side lengths fixed), see  \cite{PanGordTep}, \cite{Pankhi},\cite{panzh}, \cite{zhu}; {for the space} with perimeter fixed \cite{ropepaper}  and with other constraints. Probably the most complete case has been treated in \cite{Dan}.

\medskip
\subsection*{Polygons with fixed perimeter \cite{ropepaper},\cite{Dan}}
\begin{figure}[h]
\centering \includegraphics[width=8 cm]{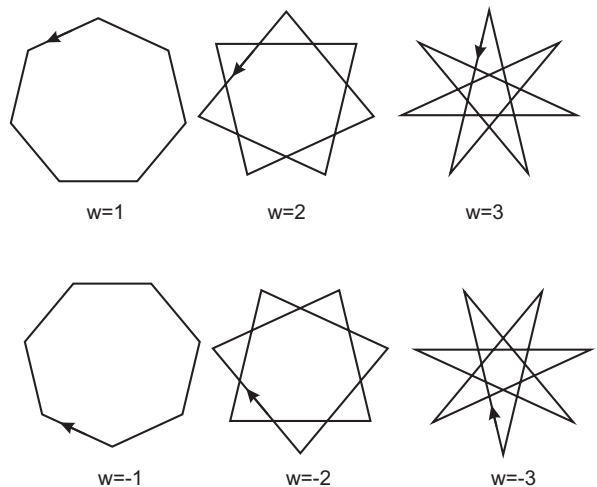}
\caption{Regular stars for $n=7$ with their winding numbers.}\label{Star7}
\end{figure}

We shall need the following result:

 Let  $\sigma$ be a cyclic renumbering: given a polygon $P=(p_1,...,p_n)$, $$\sigma(p_1,...,p_n)=(p_2,p_3,...,p_n,p_1).$$ In other words, we have an action of $\mathbb{Z}_n$ on $\mathcal{C}$  which renumbers the vertices of a polygon cyclically.
  A \textit{regular star} is an equilateral $n$-gon  such that $\sigma (p_1,...,p_n)=(p_1,...,p_n)$, see Figure \ref{Star7}.

  A \textit{complete fold}  is a regular star with $p_i=p_{i+2}$. It exists for even $n$ only.

 A regular star $P$ which is not a complete fold is uniquely defined by its winding number $w(P)$ with respect to the center.

Let $\mathcal{C}_{\pi}$ be the space of all $n$-gons with some fixed  value $\pi>0$ of perimeter.

\newpage

\begin{theorem} \cite{ropepaper},\cite{Dan} \label{CritInd} \begin{enumerate}
\item $\mathcal{C}_\pi$ is homeomorphic to $\mathbb{C}P^{n-2}$.
 \item  {The restriction} $\mathcal{A}|_{\mathcal{C}_\pi}$ is a perfect Morse function.

                   \item The critical points of the function $\mathcal{A}|_{\mathcal{C}_\pi}$ are regular stars and complete folds only.

                   \item The Morse index of a regular star is:
                   $$M(P)=\left\{
                            \begin{array}{ll}
                              2w(P)-2, & \hbox{if $w(P)<0$;} \\
                             2n-2w(P)-2, & \hbox{if $w(P)>0$;} \\
                              n-2, & \hbox{if $P$ is a complete fold (for even $n$ only).}
                            \end{array}
                          \right.
                   $$

                   \item  The critical values of function $\mathcal{A}|_{\mathcal{C}_\pi}$ increase with the Morse index.
                 \end{enumerate}
\end{theorem}
The theorem was proved in \cite{ropepaper} under assumption that the critical points are non-degenerate Morse points. Later on, it was proven in \cite{Dan} that the the critical points are indeed non-degenerate. The proof was missing for complete folds, but now we provide the proof in the Appendix.

\medskip

It follows from (1) that the space $\mathcal{C}$ of all $n$-gons   can be viewed as a cylinder over $\mathbb{C}P^{n-2}$, or rather as a cone over $\mathbb{C}P^{n-2}$ with the apex removed.

\bigskip
{\subsection*{ Polygons with fixed area}
Denote by $\mathcal{C}_a$ the space of polygons with fixed value $a$ of $\mathcal{A}$. (We distinguish between $\mathcal{C}_a$
and $\mathcal{C}_\pi$ by arabic and greek letters.) By scaling, it is clear that all $\mathcal{C}_a$ with $a\neq 0$ are homeomorphic.
In Theorems 4 and 5 we state analogues  of Theorem \ref{CritInd} for $\mathcal{P}|_{\mathcal{C}_a}$.}

\bigskip

\subsection*{ Area and perimeter in interaction}

The space $\mathcal{C}$ has two foliations: one by level sets of the perimeter, and the other one by level sets of the oriented area. Area can be considered as a function on the leaves of the perimeter foliation,  and vice versa.

One of the main goals of the paper is to compare the level foliations of $\mathcal{A}$ and $\mathcal{P}$.
{Let us first mention that the minimum of $\mathcal{P}|_{\mathcal{C}_a}$ corresponds to the maximum of $\mathcal{A}|_{\mathcal{C}_\pi}$. This is called the 'isoperimetric duality'. To extend it to all critical points and level surfaces, let us take together
 area and perimeter}.  We obtain   a map
$$\Phi:\mathcal{C}\rightarrow \mathbb{R}_{>0}\times \mathbb{R},$$

$$\Phi(P)=(\mathcal{P}(P),\mathcal{A}(P)).$$

\begin{proposition}\begin{enumerate}
                     \item The critical set $\Gamma \subset \mathcal{C}$ of $\Phi$ consists of finitely many curves.
                     \item The discriminant set $D \subset \mathbb{R}_{>0}\times \mathbb{R}$ of $\Phi$ is a finite set of curves $D_i$  given by $\mathcal{A}=c_i\mathcal{P}^2$. The curves  $D_i$ correspond to regular stars and complete folds.
                     \item $\Phi$ is a locally trivial fiber bundle over $\mathbb{R}_{>0}\times \mathbb{R}\setminus D$.
                     \item The regular fibers are submanifolds of $\mathcal{C}$ of codimension $2$. Each singular fiber has a Morse-type singularity.
                   \end{enumerate}
\end{proposition}
Proof. If a polygon $P$ is critical for $\Phi$, then so is any dilate of $P$. Due to homogeneity of $\mathcal{A}$ and $\mathcal{P}$,  the map $\Phi$ sends each line $cP$ to a curve $A=cP^2$ in $\mathbb{R}_{>0}\times \mathbb{R}$. The singularities of $\mathcal{A}$  on $\mathcal{C}_{\pi}$ are isolated and of Morse type, see Theorem \ref{CritInd}. Therefore $\Gamma$ is one-dimensional, and the fibres above points of $D$ are of Morse type. The index is defined for every transversal slice, and due to homogeneity  property, the index is constant for each curve.\qed

\medskip

\begin{figure}[h]
\centering \includegraphics[width=16 cm]{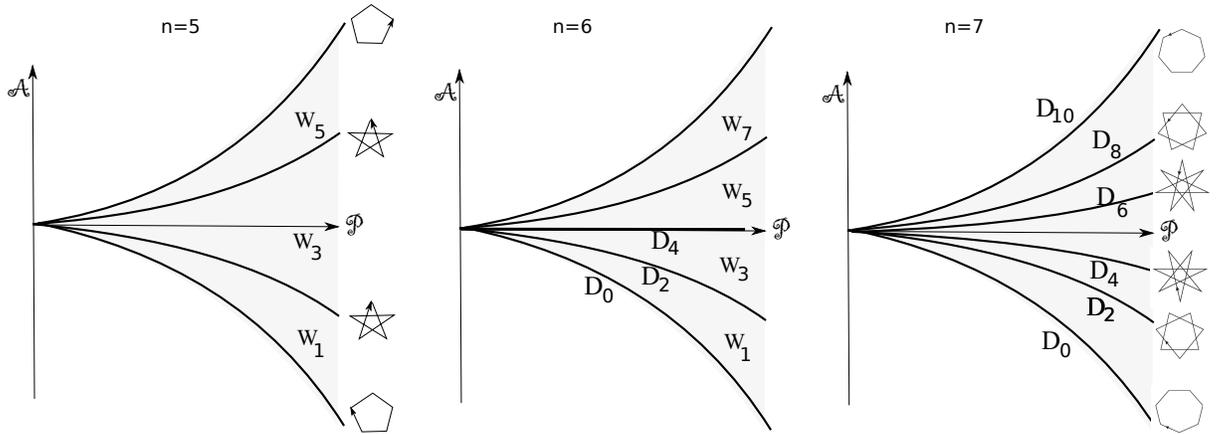}
\caption{Stratification of $\mathbb{R}_{>0}\times \mathbb{R}$.}\label{Fig2}
\end{figure}

To summarize, one can envisage the situation as follows. The image of $\Phi$ is the shadowed area in Figure \ref{Fig2}.

 The bold curves $D_i$  are the images of the critical points. We number them  accordingly to the Morse indices. For even and odd $n$ we have two different pictures: if $n$ is even, there appears a horizontal bold ray which corresponds to the complete folds. Each vertical straight line is  the image of the space $\mathcal{C}_{\pi}$ for some $\pi$,
and its intersection with a bold curve corresponds to a critical value described in Theorem \ref{CritInd}.

The bold curves split the image of $\Phi$ into open domains; let us enumerate them starting from the bottom by odd numbers: $W_1,W_3,W_5,$ etc. In each domain we have the same fiber type. Together with singular fibers, they constitute a codimension $2$ foliation over $\mathcal{C}$ which is fiber homeomorphic to a product of $\mathbb{R}_{>0}$ {(foliated by points)} and the $\mathcal{A}$-foliation on $\mathcal{C}_\pi$.

These constructions and concepts are known as polar curves and Cerf diagrams in (complex) singularity theory, see e.g. \cite{Tibar}.

\section{Configuration spaces of polygons with area and perimeter fixed.}

For a point $(\pi,a)\in \mathbb{R}_{>0}\times \mathbb{R}$ set $\mathcal{C}_{\pi,a}=\Phi^{-1}(\pi,a)$
\begin{theorem}\label{ThmOdd}
\begin{enumerate}
                 \item Let $i\leq n-2$. Assume that $(\pi,a)\in W_{i}$.
The fiber $\mathcal{C}_{\pi,a}$ over this point is an $2n-5$-dimensional non-singular manifold whose topological type depends on the domain $W_i$ only.
The homology groups of $\mathcal{C}_{\pi,a}$ are:

 $$H_{j}(\mathcal{C}_{\pi,a},\mathbb{Z})=\left\{
                                 \begin{array}{lll}
                                   \mathbb{Z}, & \hbox{if $j$ is even and   $j\leq i$ ;  } \\
                                   \mathbb{Z}, & \hbox{if $j$ is odd and   $2n-5\geq j\geq 2n-5-i $ ;  } \\
                                   0, & \hbox{otherwise.}
                                 \end{array}
                               \right.$$

\item { For $(\pi,a)\in D_0$, the space $\mathcal{C}_{\pi,a}$  is a single point. Let $0<i<n-2$, and  $(\pi,a)\in D_i$.} Then $\mathcal{C}_{\pi,a}$  is an $2n-5$-dimensional manifold with one singular point.
Its homology groups are:

$$H_{j}(\mathcal{C}_{\pi,a},\mathbb{Z})=\left\{
                                 \begin{array}{lll}
                                   \mathbb{Z}, & \hbox{if $j$ is even and   $j\leq i$ ;  } \\
                                   \mathbb{Z}, & \hbox{if $j$ is odd and   $2n-5 \geq j\geq 2n-5-i $;  } \\
                                   0, & \hbox{otherwise.}
                                 \end{array}
                               \right.$$

 \item If $i> n-2$, { then
                  the fiber  $\mathcal{C}_{\pi,a}$ of a point $(\pi,a)\in W_i$ (respectively, of a point $(\pi,a)\in D_i$) is homeomorphic to the fiber over a point from  $W_{2n-i-4}$ (respectively, from  $D_{2n-i-4}$).}

{Their homology groups follow from (1) and (2).}

\medskip

 \item    $\mathcal{C}_{\pi,0}$  is an $2n-5$-dimensional manifold (with a singularity if $n$ is even).
Its homology groups are:

 $$H_{j}(\mathcal{C}_{\pi,0},\mathbb{Z})=\left\{
                                 \begin{array}{lll}
                                   \mathbb{Z}, & \hbox{if $j$ is even and   $j\leq n-2$ ;  } \\
                                   \mathbb{Z}, & \hbox{if $j$ is odd and   $2n-5\geq j>n-2$;}\\
                                   0,& \hbox{otherwise.}
                                 \end{array}
                               \right.$$

               \end{enumerate}

\end{theorem}
Proof.

 (3) {is true since} symmetry with respect a line keeps perimeter and takes $\mathcal{A}$ to $-\mathcal{A}$.

 Let us \textbf{prove (1)}. Consider the restriction   $\mathcal{A}|_{\mathcal{C}_\pi}$. Set $A=\Phi^{-1}(\pi, (-\infty, a])$. By Morse theory, it has the homotopy type of a cell complex with even-dimensional cells only, one cell in each dimension starting from  $0$ to $i-1$.

 Let us explain this in more details. Morse theory suggest to control the behaviour of $A=\mathcal{A}|_{\mathcal{C}_\pi}^{-1}(- \infty,x]$ and
  $\mathcal{A}^{-1}(x)$
  as $x$ grows. When $x$ is very small, both preimages are empty. The space $A$ becomes a ball after passing through the minimum of  $\mathcal{A}|_{\mathcal{C}_\pi}$, whereas $\mathcal{A}|_{\mathcal{C}_\pi}^{-1}(x)$ turns to a sphere. For a regular value $x$ is a non-singular manifold. Passing through a critical value means (a) attaching a cell  to $A$ whose dimension equals the Morse index  and (b) applying a Morse surgery  to $\mathcal{A}|_{\mathcal{C}_\pi}^{-1}(x)$.
   In our case  all the indices  are even.

 An important feature of our particular case is that the function $\mathcal{A}|_{\mathcal{C}_\pi}$ is not everywhere continuously differentiable, since the perimeter involves square roots.  This problem has been treated in \cite{ropepaper}. In short,
  Clarke subdifferential  from nonsmooth analysis  generalizes smooth Morse theory
and provides  the "regular interval"  theorem which  states, that (1) preimage of a Lipshitz-regular value is a manifold,
and (2) if an interval contains no critical values, then the preimages of its endpoints are homeomorphic.

In \cite{ropepaper} we  checked that at each non-smooth point  $\mathcal{A}|_{\mathcal{C}_\pi}$ is regular, see Theorem \ref{CritInd}.

 Analogously, the set $B=\Phi^{-1}(\pi, [a,\infty))$ has the homotopy type of a cell complex with even-dimensional cells only, one cell in each dimension starting from  $0$ to $2n-i-5$.

 The union of $A$ and $B$ is homeomorphic to $\mathbb{C}P^{n-2}$ (Theorem \ref{CritInd}), and the intersection is exactly the space $\mathcal{C}_{\pi,a}$ we are interested in.  The Mayer-Vietoris exact sequence for these spaces
  can be visualized as a table, where the rows correspond to the indices of homology groups, and whose entries are
 homology groups of the participants of the sequence. The table below depicts the case $n=7,$ and $(\pi,a)\in W_3$.
 Initially, the first column contains unknown homology groups. Let us show that they are recovered uniquely.

\begin{lemma}\label{LemmaId} Let $\mathbf{i}^A: A\subset A\cup B$ and  $\mathbf{i}^B: B\subset A\cup B$ be inclusions. Then the induced  homomorphisms $$\mathbf{i}^A_*:H_*(A, \mathbb{Z})\rightarrow H_*(A\cup B, \mathbb{Z})\hbox{ \ and  \ } \mathbf{i}^B_*:H_*(B, \mathbb{Z})\rightarrow H_*(A\cup B, \mathbb{Z}) $$ are  monomorphisms.
\end{lemma}
Proof follows from the cell complex structures of $A$, $B$, and $A\cup B$: up to dimension $i$, the complex $A\cup B$
has exactly the same cells as $A$.

  \bigskip

  Depending on the index $i$, the  Mayer-Vietoris sequence splits into three types of sequences:

 (1) In the upper part, for $j=2k, \ \  2n-5-i\leq j\leq 2n-4$,  we have
$$0\rightarrow H_{2k}(A\cap B)\rightarrow 0 \rightarrow \mathbb{Z}\rightarrow H_{2k-1}(A\cap B)\rightarrow 0.$$
One concludes that  $H_{2k}(A\cap B)=0,\ \ H_{2k-1}(A\cap B)=\mathbb{Z}. $

\medskip

  (2) Next comes $j=2k, \ \ i\leq j\leq 2n-5-i$.

  $$0\rightarrow H_{2k}(A\cap B)\rightarrow \mathbb{Z} \rightarrow \mathbb{Z}\rightarrow H_{2k-1}(A\cap B)\rightarrow 0$$
{By Lemma \ref{LemmaId}} the central arrow is the identity map, so the two unknown  homology groups are zero.

\medskip

(3) Finally,  for $j\leq i$, we have
 $$0\rightarrow H_{2k}(A\cap B)\rightarrow \mathbb{Z}\oplus \mathbb{Z}  \rightarrow \mathbb{Z}\rightarrow H_{2k-1}(A\cap B)\rightarrow 0.$$

By Lemma \ref{LemmaId}, the central arrow is a surjective mapping, so $$H_{2i}(A\cap B)=\mathbb{Z},\ \ H_{2i-1}(A\cap B)=0. $$

Note that since $\mathcal{C}_{\pi,a}=A\cap B$ is a manifold, one can use the Poincar\'{e} duality and   needs not  analyze the third type of the sequences.

\bigskip

{Example: the table} for $n=7,$ and $(\pi,a)\in W_3$ is completed like this:

\medskip

\small{ \begin{tabular}{ | l | l | l | l | }
\hline
j&$H_*(\mathcal{C}_{\pi,a},\mathbb{Z})$ & $H_*(A,\mathbb{Z})\oplus H_*(B,\mathbb{Z})$ & $H_*(A\cup B,\mathbb{Z})$ \\ \hline

10 & 0 & $0$ & $\mathbb{Z}$ \\
   9 & $\mathbb{Z}$  & 0 & 0 \\
   8 & 0  & $0 $ & $\mathbb{Z}$ \\
   7 & $\mathbb{Z}$  & 0 & 0 \\
   6 & 0  & $0 \oplus \mathbb{Z}$ & $\mathbb{Z}$ \\
   5 & 0  & 0 & 0 \\
   4 & 0  & $ 0 \oplus \mathbb{Z}$ & $\mathbb{Z}$ \\
   3 & 0 & 0 & 0 \\
   2 & $\mathbb{Z}$ &  $\mathbb{Z}\oplus \mathbb{Z}$ & $\mathbb{Z}$ \\
   1 & 0 & 0 & 0 \\
   0 & $\mathbb{Z}$ & $\mathbb{Z}\oplus \mathbb{Z}$ & $\mathbb{Z}$ \\

\hline
\end{tabular}}

\bigskip

 Now let us \textbf{prove (2) and (4)}. Set $A=\Phi^{-1}(\pi, (-\infty, a+\varepsilon])$ with some small $\varepsilon >0$. By Morse theory, $A$ has a homotopy type of a cell complex with even-dimensional cells only, one cell in each even dimension ranging from  $0$ to $i$.
 Set also $B=\Phi^{-1}(\pi, [a-\varepsilon,\infty))$. It also has one cell in each even dimension ranging from $0$ to $2n-4-i$.

 The intersection $A\cap B$ is homotopy equivalent to $\mathcal{C}_{\pi,a}$, and a similar analysis of
the Mayer-Vietoris sequence completes the proof.

\bigskip

Example: the table below depicts the case $n=6, \ a=0, \ i=4$.

\medskip

 \begin{tabular}{ | l | l | l | l | }
\hline
&$H(\mathcal{C}_{\pi,0},\mathbb{Z})$ & $H(A,\mathbb{Z})\oplus H(B,\mathbb{Z})$ & $H(A\cup B,\mathbb{Z})$ \\ \hline

   8 & 0  & $0 $ & $\mathbb{Z}$ \\
   7 & $\mathbb{Z}$  & 0 & 0 \\
   6 & 0  & $0$ & $\mathbb{Z}$ \\
   5 & $\mathbb{Z}$  & 0 & 0 \\
   4 & $\mathbb{Z}$  & $ \mathbb{Z} \oplus \mathbb{Z}$ & $\mathbb{Z}$ \\
   3 & 0 & 0 & 0 \\
   2 & $\mathbb{Z}$ &  $\mathbb{Z}\oplus \mathbb{Z}$ & $\mathbb{Z}$ \\
   1 & 0 & 0 & 0 \\
   0 & $\mathbb{Z}$ & $\mathbb{Z}\oplus \mathbb{Z}$ & $\mathbb{Z}$ \\

\hline
\end{tabular}

\qed

\bigskip

\section{Duality in the isoperimetric problem}

The following theorem demonstrates equivalence of level foliations of $\mathcal{A}$ and $\mathcal{P}$.

\begin{theorem}\label{ThmRotate} If $a> 0$, there exists a level-preserving homeomorphisms
$$\Psi_+:\mathcal{C}_{a}\longrightarrow  \mathcal{C}_{\pi}|_{\mathcal{A}>0}$$ which sends Morse points of index $i$ to those of index $2n-4-i$.

If $a<0$, the analogous homeomorphism 
$$\Psi_-:\mathcal{C}_{a}\longrightarrow  \mathcal{C}_{\pi}|_{\mathcal{A}<0}$$  preserves Morse indices.
\end{theorem}

Proof.  {Assume that $a>0$.  We need to construct the horizontal maps in the following commutative diagram:}
$$\begin{array}{ccc}
  \mathcal{C}_a & \cong & \mathcal{C}_{\pi}|_{\mathcal{A}>0} \\

  \downarrow \mathcal{P} &  & \downarrow \mathcal{A} \\
  \mathbb{R}_{>0} & \cong & \mathbb{R}_{>0}.
\end{array}
$$

\begin{figure}[h]
\centering \includegraphics[width=10 cm]{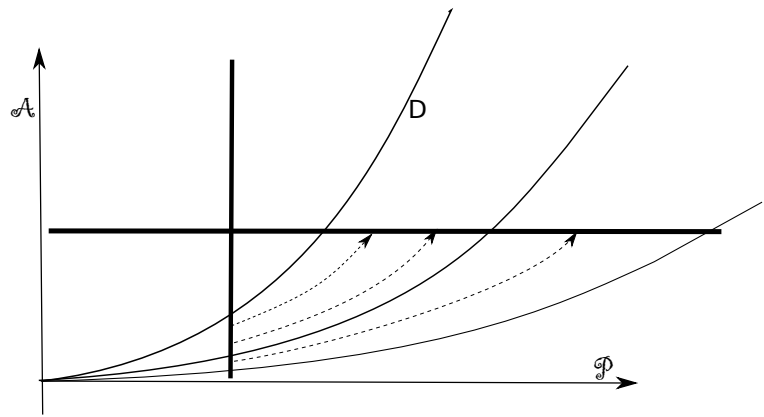}
\caption{Rotating the  half line  to the horizontal position}\label{FigRotate}
\end{figure}

  Observe that $\Phi$ is a stratified submersion. The proof comes from "le rabattement de diagramme de Cerf", well known in singularity theory, and used e.g. in \cite{Tibar}  and \cite{DirkMonodr}.
We rotate the vertical half line $\mathcal{P}=\pi$ with $\mathcal{A}>0$ to the horizontal position $\mathcal{A}=a$.
An explicit map $$\mathcal{C}_\pi\cap \{\mathcal{A}>0\}\rightarrow  \mathcal{C}_a$$  can be obtained by scaling action, which gives a level preserving homeomorphism, see Figure \ref{FigRotate}. It follows that $\mathcal{P}$ and $-\mathcal{A}$ have the same Morse evolution.

{Since $\mathcal{A}=0$ splits $\mathcal{C}_\pi$ into two symmetric pieces, the analogous statement applies for $\mathcal{C}_{\pi}|_{\mathcal{A}<0}$.}
\qed

\bigskip

The proof of Theorem \ref{ThmRotate} reveals the following general statement about 'kissing' surfaces:

\begin{proposition}
Let $$f:(\mathbb{R}^{n+1},0)\rightarrow (\mathbb{R},0) \hbox{\ \ and \ \ } g:(\mathbb{R}^{n+1},0)\rightarrow (\mathbb{R},0) $$
be germs of two functions such that $X=f^{-1}(0)$ and $Y=g^{-1}(0)$ are smooth hypersurfaces tangent to each other at $(0,0)$.
Assume that $f|_Y$ has a Morse singularity of index $M$ at $(0,0)$. Then also $g|_X$ has a Morse singularity at $(0,0)$, and its index $m$ satisfies
$$
\left\{
  \begin{array}{ll}
  M+m=n,  & \hbox{if \  $\nabla f$ and $\nabla g$ are codirected;} \\
   M=m, & \hbox{if \  $\nabla f$ and $\nabla g$ are oppositely directed.}\qed
  \end{array}
\right.
$$
\end{proposition}

\medskip

NB. This fact can be also demonstrated by analytic methods. It serves as the first important step in understanding the duality.

\medskip

We immediately obtain the "dual" version of Theorem \ref{ThmOdd}:

\begin{theorem}\label{ThmDual}For $a\neq 0$, the restriction $\mathcal{P}|_{\mathcal{C}_a}$  is a perfect Morse function.
\begin{enumerate}
                 \item If $a>0$, the critical points of $\mathcal{P}|_{\mathcal{C}_a}$ are exactly the regular stars with positive winding numbers.
For each of the critical points, its Morse index
satisfies { $$m(P)+ M(P)=2n-4.$$}

                 \item If $a<0$, the critical points  of $\mathcal{P}|_{\mathcal{C}_a}$are exactly the regular stars with negative winding numbers.
For each of the critical points, its Morse index
satisfies $$m(P)= M(P).$$
{In both cases  $M(P)$ comes from Theorem \ref{CritInd}.}
               \end{enumerate}

\end{theorem}

Another theorem follows directly.

\begin{theorem} \label{ThmAreaFixed}\begin{enumerate}
                          \item The space $\mathcal{C}_a$ of polygons with prescribed area $a\neq 0$ is a $2n-4$-dimensional manifold and has a homotopy type of a cell complex with even cells only. The dimensions of the cells range from $0$ to $[\frac{n+1}{2}]\cdot2-4$, one cell in each even dimension.

                              \item The space of polygons with zero area  is homeomorphic to $\mathcal{C}_{\pi,0}\times \mathbb{R}_{>0}$, so its homologies follow from Theorem \ref{ThmOdd}.\qed

                    \end{enumerate}
\end{theorem}

\bigskip

Let us finish the paper by the \textbf{concluding remark:}

   A similar  duality approach is possible for any two functions $f$ and $g$, especially if both are homogeneous. In particular, one can think about sum of squared edge lengths versus area on the polygon space.

Another example is
Coulomb energy and perimeter on polygon space (in this case critical points for $n=3$ and fixed perimeter are described in
\cite{KPSThree}. An interesting question is to extend the approach of the present paper to three or more functions.

\section{Appendix}

\begin{proposition}
 The complete fold is a Morse point of the function $\mathcal{A}$ restricted to $\mathcal{C}_\pi$.
\end{proposition}
We may assume that the complete fold is $P^o=(p^o_1,...,p^o_n)$  with $p_1^o=0$ such that
$\overrightarrow{a}_i^o= \overrightarrow{p^o_ip^o}_{i+1}=(-1)^{i+1}\overrightarrow{a}$ for some non-zero vector $\overrightarrow{a}$.
Let us coordinatize a neighborhood of $P^o$  by $(\overrightarrow{b}_1,...,\overrightarrow{b}_{n-1})$ where

$\overrightarrow{a}_i= \overrightarrow{p_ip}_{i+1}=(-1)^{i+1}\overrightarrow{a}+ \overrightarrow{b}_i, \ \ \ i=1,...,n-1.$

It is important that these coordinates describe the space of polygons with rotations and dilations not factored out.

Since the area vanishes at $P^o$, computing the $2$-jet of $\frac{\mathcal{A}}{\mathcal{P}^2}$ reduces to computing the $2$-jet of $\mathcal{A}$.

Now direct computations show that the  Hessian matrix of $\frac{\mathcal{A}}{\mathcal{P}^2}$ has the following form:

$$
\left(
  \begin{array}{cccccccccc}

    0 & 0 & 0 & 1 & 0 & 0 & 0 & 0 & 0 & 0 \\
    0 & 0 & -1 & 0 & 0 & 0 & 0 & 0 & 0 & 0 \\
    0 & -1 & 0 & 0 & 0 & 1 & 0 & 0 & 0 & 0 \\
    1 & 0 & 0 & 0 & -1 & 0 & 0 & 0 & 0 & 0 \\
    0 & 0 & 0 & -1 & 0 & 0 & 0 & 1 & 0 & 0 \\
    0 & 0 & 1 & 0 & 0 & 0 & -1 & 0 & 0 & 0 \\
    0 & 0 & 0 & 0 & 0 & -1 & 0 & 0 & 0 & 1 \\
    0 & 0 & 0 & 0 & 1 & 0 & 0 & 0 & -1 & 0 \\
    0 & 0 & 0 & 0 & 0 & 0 & 0 & -1 & 0 & 0 \\
    0 & 0 & 0 & 0 & 0 & 0 & 1 & 0 & 0 & 0 \\

  \end{array}
\right)
$$

This is the matrix for $n=6$. In the general case we have a similar matrix $2(n-1)\times 2(n-1)$.
 It is sufficient to prove that its rank equals $2n-4$ since two zero eigenvalues correspond to rotations and dilations.
Eliminate the last two rows and the last to columns of the matrix. The resulted matrix is non-degenerate. Indeed, it is easy to conclude that its kernel vanishes.\qed

\section*{Acknowledgments} This work is
 supported by the RFBR grant 20-01-00070.
It is our pleasure to acknowledge the  hospitality and excellent working conditions of
CIRM, Luminy, where
this paper was initiated and completed  as a 'research in pairs' project.

\end{document}